\documentclass[12pt]{amsart}
\newtheorem{thm}{Th\'eor\`eme}[section]% TH\'EOR\`EME

\newtheorem{lem}[thm]{Lemme}% LEMME

\newtheorem{remark}{Remarque} % REMARQUE

\begin{document}
{
\center{\Large{\bf 
Syst\`eme des coordonn\'ees plates de la 
d\'eformation 
verselle des singularit\'es isol\'ees simples d'hypersurface
et z\'eros de l'int\'egrale d'Abel 

} }
 \vspace{1pc}
\center{\large{ Susumu TANAB\'E}}
 }
 \vspace{0.5pc}

\begin{center}
\begin{minipage}[t]{12.2cm}
{\sc R\'esum\'e -} {\em
On \'etudie le syst\`eme des \'equations differ\'entielles satisfaites
par l'int\'egrale d'Abel associ\'ee au cycle
 $\gamma_s \subset \{(x,y) \in {\bf R}^2:
H(x,y;s)
=0 \}$ d\'efinie pour la d\'eformation verselle d'une singularit\'e isol\'ee 
simple d'hypersurface.   
Comme application,
on obtient une estimation de la multiplicit\'e  des z\'eros  de
l'int\'egrale $I_{\omega}(s)=
\int_{\gamma_s} \omega $ en fonction de poids quasihomog\`enes associ\'es
\`a $H(x,y;0)$ et de $deg (\omega).$
 }
\end{minipage}
\end{center}

%%%%%%%%%%%%%%%%%%%%%%%%%%%%% SECTION 1%%%%%%%%%%%%%%%%%%%%%%%%%%%%%%
{\center{\section{\bf 
Introduction
}}}

Dans cette note, nous poursuivons le probl\`eme de l'estimation des
 z\'eros de l'int\'egrale d'Abel \'etudi\'e par plusieurs auteurs ~\cite{Gav2}, ~\cite{Mar} 
etc. C'est une demarche vers une reponse raisonnable au XVIe probl\`eme 
de Hilbert sur le nombre des cycles limites. Dans un travail pr\'ecedent
~\cite{T98}, nous avons \'etabli une estimation de la 
multiplicit\'e des z\'eros de l'int\'egrale hyperelliptique associ\'ee \`a la singularit\'e $A_\mu,$ en se servant du syst\`eme de Gau\ss-Manin satisfait par l'int\'egrale hyperelliptique. La strat\'egie de cette note est semblable 
\`a celle de ~\cite{T98}, pourtant la tactique en diff\`ere au sens que nous utilisons ici l'expression du syst\`eme de Gau\ss-Manin au moyen des coordonn\'ees plates introduites par ~\cite{SSY}. 
Dans ~\cite{SSY}, les coordonn\'ees plates (flat coordinates) ont \'et\'e introduites comme les invariants canoniques d'un groupe de Coxeter fini 
irreductible qui satisfont certaines conditions naturelles. C'est gr\^ace \`a Noumi ~\cite{No} que l'utilit\'e pr\'epond\'erante de ces coordonn\'ees se manifeste aux \'etudes du syst\`eme de Gau\ss-Manin. Notre demarche s'appuie sur un r\'esultat de ~\cite{No} qui est formul\'e au Th\'eor\`eme ~\ref{thm21}. Ce 
r\'esultat nous permet d'exprimer avec les coordonn\'ees plates, le syst\`eme de de Gau\ss-Manin associ\'e \`a la d\'eformation verselle des singularit\'es isol\'ees simples d'hypersurfaces, d'une fa\c{c}on plus simple que l'expression 
du syst\`eme par les 
coordonn\'ees usuelles de l'espace de d\'eformation verselle 
(voir ~\cite{T98}).   
 L'objet central de notre recherche est l'int\'egrale d'Abel 
associ\'ee \`a  une singularit\'e isol\'ee simple dhypersurface pour
laquelle on \'etablit le r\'esultat suivant:

{\bf Th\'eor\`eme}
{\it Soit
$\gamma_s \subset \{(x,y) \in \mathbf R^2 : H(x,y;s)=0\}$
un cycle \'evanescent de la d\'eformation verselle d'une singularit\'e isol\'ee
simple d'hypersurface de la liste (2.1) i.e. $A_\mu,$ $D_\mu,$  $E_6,$ 
$E_7,$ $E_8$. 
Alors il existe un ensemble ouvert dense
$U \subset \mathbf R^{\mu} \setminus D$
tel que la multiplicit\'e des 
z\'eros $N$ de l'int\'egrale d'Abel (1.1) non-nulle
$$I_{P,Q,\gamma_s}(s) =\int_{\gamma_s} P(x,y)dx + Q(x,y)dy\leqno{(1.1)}$$
avec $P(x,y), Q(x,y)$ des polyn\^omes de degr\'e au plus $K$
i.e.$P(x,y) = \sum_{0 \leq i+j \leq K} P_{i,j}x^iy^j,$
$Q(x,y)$ $ = $ $\sum_{0 \leq i+j \leq K}$ $ Q_{i,j}$ $x^iy^j.$
doit satisfaire \`a chaque point
$s \in U$
l'in\'egalit\'e suivante:
$$N \le \mu (1 +  [(\rho_1 +\rho_2) (K-1)])-1, \leqno(1.2)$$
o\`u $\rho_1, \rho_2$ sont les poids quasihomog\`enes des variables 
$x,y$ d\'efinis pour chaque singularit\'e.
}

Il est facile de d\'eduire de $(1.2)$ que si 
$\mu,\nu, K ,m \leq n, $
alors la multiplicit\'e des z\'eros $N$ est domin\'e par une fonction 
quadratique en  $n.$ 
Nos r\'esultats, donc,  paraissent une  am\'elioration  des
estimations obtenues dans  ~\cite{Mar} ( la multiplicit\'e de chaque
z\'ero $ \leq \frac{n^4 + n^2 -2}{2}$  )  
pour autant qu'il s'agisse de l'int\'egrale du type  $(1.1).$  

\vspace{0.5pc}
%%%%%%%%%%%%%%%%%%%%%%%%%%%%% SECTION 2%%%%%%%%%%%%%%%%%%%%%%%%%%%%%%
  {\section {\bf D\'eformation verselle 
des singularit\'es isol\'ees simples d'hypersurface} }

D'abord, on regarde la d\'eformation verselle des 
singularit\'es isol\'ees simples d'hypersurface:
$$\begin{array}{l}
A_{\mu} : H(x,y;s)  =  x^{\mu+1} + y^2 +\sum^{\mu-1}_{i=0} s_{i0} x^i, \\
D_{\mu} : H(x,y;s)  =  x^{\mu-1} + xy^2 + \sum^{\mu-2}_{i=1} s_{io}x^{i-1} + s_{01}y + s_{00}, \\
E_6 : H(x,y;s)  =  x^4 + y^3 + s_{21}x^2y+s_{11}xy + s_{20}x^2 + s_{10}x + s_{01}y + s_{00}, \\
E_7 : H(x,y;s)  =  x^3 + xy^3 + s_{21} x^2y + s_{11}xy + s_{20}x^2 + s_{02} y^2 + s_{10}x + s_{01}y + s_{00}, \\
E_8 : H(x,y;s)  =  x^5+y^3 + s_{31}x^3y+s_{21}x^2y+s_{11} xy+s_{01}y + s_{30}x^3 + s_{20}x^2 + s_{10}x + s_{00}.
\end{array} \leqno{(2.1)}$$
Nous d\'efinissons un ensemble fini
$\mathbf I \subset \mathbf N^2$
associ\'e \`a chaque singularit\'e  ci-dessus qui satisfait  les 
conditions suivantes:

\begin{enumerate}
\item $\mid \mathbf I \mid = \mu =$ le nombre de Milnor de la singularit\'e.
\item A chaque \'el\'ement $(\nu_1,\nu_2) \in \mathbf I$ correspond un mon\^ome $x^{\nu_1} y^{\nu_2}$ qui est un \'el\'ement de la base de l'espace

$$\mathbf C[x,y] /(\frac {\partial} {\partial x} H(x,y;0)),
 \frac {\partial} {\partial y} H(x,y;0)).$$
\end{enumerate}

On remarque que toutes les singularit\'es de la liste (2.1) sont quasihomog\`enes. C'est \`a dire, il existe
$(\rho_1, \rho_2) \in \mathbf Q^2_{>0} $ tel que

$$(\rho_1x \frac{\partial} {\partial x} + \rho_2y \frac{\partial} {\partial y}) 
H(x,y;0) = H(x,y;0).$$
On d\'efinit pour
$\alpha \in \mathbf N^{\mu}$
le poids  quasihomog\`enes d'un mon\^ome
$s^{\alpha}$:
$$weight \left( \prod_{(\nu_1,\nu_2) \in \mathbf I} s^{\alpha_{\nu_1,\nu_2}} \right) = \langle w,\alpha \rangle = \sum_{(\nu_1,\nu_2) \in \mathbf I} w_{\nu_1,\nu_2} \alpha_{\nu_1\nu_2} \leqno{(2.2)}$$
o\`u
$w_{\nu_1\nu_2} = 1-(\rho_1\nu_1+\rho_2\nu_2)$.
On d\'efinit pour chaque
$(\nu_1,\nu_2) \in \mathbf I$
un r\'eseau d'entiers
$L(\nu_1,\nu_2)$
selon Noumi \cite{No}. Dans les cas de
$A_{\mu}, E_6, E_8:$

$$L(\nu_1,\nu_2) = ((\nu_1,\nu_2) + L) \cap \mathbf N^2$$
o\`u
$L = \mathbf Z (\frac 1 {\rho_1}, 0) + \mathbf Z (0, \frac 1 {\rho_1}) \subset \mathbf Z^2$.
Soit
$$
C_{\nu_1\nu_2}(\beta_1,\beta_2)  =  (-1)^{\rho_1(\beta_1-\nu_1)+\rho_2(\beta_2-\nu_2)} \frac {\Gamma(\rho_1(\beta_1+1))\Gamma(\rho_2(\beta_2+1))} {\Gamma(\rho_1(\nu_1+1))\Gamma (\rho_2(\nu_2+1))}\qquad\mbox{si}\qquad (\beta_1,\beta_2) \in L(\nu_1,n_2) 
\leqno{(2.3)}$$
$$=  0 \qquad\mbox{si}\qquad (\beta_1,\beta_2) \notin L(\nu_1,\nu_2). $$
Dans les cas de
$D_{\mu}, E_7 : L = \mathbf Z(\frac 1 {\rho_1}, 0)+ \mathbf Z(1,\frac 1 {\rho_2}),$
$L(\nu_1,\nu_2) = ((\nu_1,\nu_2)+L) \cap \mathbf N^2$.
Soit
\begin{flushleft}
$
(2.4) \;\;C_{\nu_1\nu_2}(\beta_1,\beta_2)$ 
$$=  (-1)^{\rho_1(\beta_1-\nu_1)+\rho_2(1-\rho_1)(\beta_2-\nu_2)} 
\frac {\Gamma(\rho_1[\beta_1+1-\rho_2(\beta_2+1)])
\Gamma(\rho_2(\beta_2+1))} {\Gamma(\rho_1[\nu_1+1-\rho_2(\nu_2+1)])
\Gamma (\rho_2(\nu_2+1))}
\qquad\mbox{si}\;(\beta_1,\beta_2) \in L(\nu_1,\nu_2) $$
$=  0 \qquad\mbox{si}\qquad (\beta_1,\beta_2) \notin L(\nu_1,\nu_2).$
\end{flushleft}
\par
Alors on peut introduire un syst\`eme des ``coordonn\'ees plates'' 
dans l'espace de param\`etres de d\'eformation
$\mathbf C^{\mu}_s$
(ou bien
$\mathbf R^{\mu}$) \cite{No}.
Notamment,

$$t_0 := t_{00} = s_{00} + \sum_{\langle w,\alpha \rangle =1} C_{00}(\ell (\alpha)) \frac{s^{\alpha}} {\alpha !} \leqno{(2.5)}$$

$$t_{\nu_1,\nu_2} = \sum_{\langle w,\alpha \rangle = 1-(\rho_1 \nu_1 + \rho_2\nu_2)} C_{\nu_1\nu_2} (\ell (\alpha)) \frac{s^{\alpha}} {\alpha !} \qquad\mbox{pour}\qquad (\nu_1,\nu_2) \ne (0,0) \leqno{(2.6)}$$
o\`u la fonction lin\'eaire
$\ell (\alpha)$
est d\'efinie comme suit:
$$\ell (\alpha) = \left( \sum^{\mu-1}_{i=0} i \, \alpha_{i0}, \nu_2 \right) \qquad\mbox{pour}\qquad A_{\mu}$$
$$\ell (\alpha) = \left( \sum_{(\nu_1,\nu_2) \ne (0,0)} \nu_1\alpha_{\nu_1\nu_2}, \sum_{(\nu_1,\nu_2)\ne (0,0)} \nu_2 \alpha_{\nu_1,\nu_2} \right) \qquad\mbox{pour}\qquad D_{\mu}, E_{\mu}.$$
Le changement des coordon\'ees
$s \rightarrow t$
est un diff\'eomorphisme. Desormais nous nous serverons  de la notation
$t = (t_0,t^{\prime})$,
o\`u
$t^{\prime} = (t_{\nu_1\nu_2}), (\nu_1,\nu_2) \in \mathbf I \setminus \{0,0\}$.

Soit $H(x,y;s)$
un des polyn\^omes de la liste (2.1). On prend un cycle \'evanescent
$\gamma_s \subset \{ (x,y) \in \mathbf R^2: H(x,y;s)= 0\}$.
Alors on peut d\'efinir l'int\'egrale d'Abel y associ\'ee :
$$I_{x^{k_1} y^{k_2},\gamma_s}(s) = \int_{\gamma_s} x^{k_1} y^{k_2} \frac{dx\wedge dy} {dF},$$
et un vecteur
$$\mathbf K(s) = (I_{x^{k_1}y^{k_2}}(s))_{(k_1,k_2) \in \mathbf I}.$$
L'\'enonc\'e suivant  met \`a jour l'avantage d'\'ecrire le syst\`eme
de Gauss-Manin au moyen des coordonn\'ees comme (2.5), (2.6) dans nos \'etudes.

\begin{thm}

Si on d\'efinit une base des int\'egrales d'Abel
$\mathbf J(t) = (J_{\nu_1\nu_2}(t))_{(\nu_1,\nu_2) \in \mathbf I}$
comme suit:
$$J_{\nu_1\nu_2}(t) = \int_{\gamma_{s(t)}} \frac{\partial}{\partial t_{\nu_1\nu_2}} F(x,y,s(t)) \frac{dx \wedge dy} {dF},$$
alors les \'enonc\'es suivant ont lieu.

1. Il existe une matrice
$P(t^{\prime}) \in GL (\mu, \mathbf C[t^{\prime}])$
telle que:
$$\mathbf K (s(t)) = P(t^{\prime}) \cdot \mathbf J(t) \leqno{(2.7)} $$
o\`u $P(t^{\prime}) = \frac{\partial S_{\nu}} {\partial t_{\tilde \nu}}$ 
et les indices $\tilde{\nu}, \nu \in \mathbf I$.

2. Les int\'egrales satisfont un syst\`eme d'\'equations diff\'erentielles:
$$\tilde S(t) \frac{\partial} {\partial t_0} \mathbf J(t) = \Lambda 
\cdot \mathbf J(t) \leqno(2.8)$$
o\`u
$\Lambda = \text{diag} (w_{\nu_1\nu_2})_{(\nu_1,\nu_2) \in \mathbf I}$.
La matrice
$\tilde S (t) \in \text{End} (\mathbf C^{\mu}) \otimes \mathbf C [t]$
satisfait les conditions ci-dessous:
$$\frac{\partial} {\partial t_0} \tilde S(t) = id_{\mu},$$
$$\det \tilde S (t(s)) = \Delta_{\mu}(s) \leqno{(2.9)}$$
o\`u $\Delta_{\mu}(s)$
est le discriminant de la d\'eformation verselle.
\label{thm21}
\end{thm}

On note
$D := \{s \in \mathbf C^{\mu} : \Delta_{\mu}(s) = 0\}$.
En reprenant les notations depuis le d\'ebut, on formule l'\'enonc\'e 
central  comme suit.

\begin{thm}
Soit
$\gamma_s \subset \{(x,y) \in \mathbf R^2 : H(x,y;s)=0\}$
un cycle \'evanescent de la d\'eformation verselle de la liste (2.1). Alors il existe un ensemble ouvert dense
$U \subset \mathbf R^{\mu} \setminus D$
tel que la multiplicit\'e des z\'eros $N$ de l'int\'egrale d'Abel (1.1) 
non-nulle
$$I_{P,Q,\gamma_s}(s) :=\int_{\gamma_s} P(x,y)dx + Q(x,y)dy$$
doit satisfaire \`a chaque point
$s \in U$
l'in\'egalit\'e suivante:
$$N \le \mu (1 +  [\rho_1 (K-1)+\rho_2 (K-1)])-1.$$
\label{thm22}
\end{thm}
\par
Desormais on va s'int\'eresser \`a la d\'eriv\'ee de l'int\'egrale (1.1) 
$$\frac{\partial} {\partial s_0} I_{P,Q,\gamma_s}(s) = \int_{\gamma_s}
\frac{d( P(x,y)dx + Q(x,y)dy)}{dF}
=  \int_{\gamma_s} \frac{(-P_y(x,y)+Q_x(x,y)) dx \wedge dy}{dF},$$
que l'on note par 
$$I_{R,\gamma_s}(s) 
:= \sum_{0 \leq k_1 +k_2 \leq K-1} 
\int_{\gamma_s} R_{k_1k_2} x^{k_1}y^{k_2} \frac{dx \wedge dy}{dF}.$$
On montre que la multiplicit\'e des z\'eros de l'int\'egrale 
$I_{R,\gamma_s}(s)$
ne d\'epasse pas  $\mu (1 + [\rho_1(K-1) +\rho_2 (K-1)])-2.$

Il est facile de voir qu'il existe des p\^olyn\^omes 
$P^{(k)}_{\nu_1,\nu_2}(\bullet )$ de degr\'e au plus
$v_0 := [\rho_1( K-1)+ \rho_2 (K-1)]$
pour un vecteur quelconque
$(R_{k_1k_2}) \in \mathbf R^{\frac{(K+1)K}{2}} $
(\cite {Gav2}, \cite{Yakov}):
$$\sum_{0 \leq k_1 +k_2 \leq K-1} R_{k_1k_2} I_{x^{k_1}y^{k_2}}(s(t)) = \sum_{(\nu_1,\nu_2) \in \mathbf I} P^{(k)}_{\nu_1,\nu_2}(t_0-\tilde t_0) \cdot J_{\nu_1\nu_2}(t)$$
avec 
$\tilde t_0 \in \mathbf C^{\mu} \setminus D$.

Notons le vecteur de taille
$\mu (v_0+1):$
$$\tilde{\mathbf J}(t) = ^t(J_{00}(t), \ldots, J_{\mu_1\mu_2}(t), (t_0-\tilde t_0) J_{00}(t), \ldots, (t_0-\tilde t_0) J_{\mu_1\mu_2}(t), \ldots$$
$$\ldots, (t_0-\tilde t_0)^{v_0} J_{00}(t), \ldots, (t_0-\tilde t_0)^{v_0} J_{\mu_1\mu_2} (t)).$$  
Ici
$(\mu_1,\mu_2) \in \mathbf I$
tel que
$\mu_1+\mu_2 = \max_{(\nu_1,\nu_2)\in \mathbf I} \nu_1+\nu_2$.
On note
$M = \mu(v_0+1)$.
Notre th\'eor\`eme ~\ref{thm22} se d\'eduit du lemme suivant.
\begin{lem} Il existe des polyn\^omes
$\delta^{(i)}(t), i \in \mathbf N$
de degr\'e $ \frac{M(M+1)}{2} (\mu !)^2$ tels que pour
$(\tilde t_0, t^{\prime}) \in \mathbf R^{\mu} \setminus D,$
qui se trouve hors des z\'eros de
$\delta^{(i)}(t)$,
l'\'enonc\'e suivant soit valide. 
\par
(i)Si un vecteur
$\vec{r} \in \mathbf R^M \setminus \{0\}$
satisfait la relation suivante,
$$\langle \vec{r}, \frac 1 {\ell !} \left( \frac d {dt_0}\right)^{\ell} \tilde{\mathbf J}(\tilde t_0, t^{\prime})\rangle = 0,\leqno{(2.10)}$$
pour
$0 \le \ell \le M-1$
alors (2.10) est valable pour
$\ell \ge M = \mu(v_0+1)$
aussi. 
\par
(ii) L'ensemble $\{ t \in {\bf R}^{\mu}; \delta^{(i)}(t) \not =0, i \in {\bf N} \}$
est un ouvert dense de ${\bf R}^{\mu}\setminus D.$
\label{lem23}
\end{lem}

{\bf D\'emonstration}
D'abord on \'etablit la relation entre
$\tilde{\mathbf J} (\tilde t_0, t^{\prime})$
et
$\left( \frac{\partial} {\partial t_0}\right)^{\ell} \tilde{\mathbf J} (\tilde t_0,t^{\prime})$.
Dans ce but, on introduit
$\mathbf J(t) = ^t(J_{00}(t), \ldots, J_{\mu_1\mu_2}(t))$.	
Notons le poids quasihomog\`ene d'une forme
$x^{\nu_1}y^{\nu_2}dxdy /dF$
par
$\lambda_{\nu_1\nu_2} = \rho_1(\nu_1+1)+\rho_2(\nu_2+1)-1$.
D'apr\`es (2.8) les coefficients de d\'eveloppement de Taylor de
$\mathbf J(t)$
s'\'ecrivent sous la forme suivante:
$$\left( \frac{\partial} {\partial t_0} \right)^{\ell} \mathbf J 
(\tilde t_0, t^{\prime}) = \tilde S^{-\ell}(\tilde t_0,t^{\prime}) (\Lambda -(\ell-1)id_{\mu}) \ldots (\Lambda-id) \Lambda \mathbf \cdot {\mathbf J }(\tilde t_0, t^{\prime}) \leqno{(2.11)}$$
o\`u
$\Lambda = \textnormal{diag} (\lambda_{00}, \ldots, \lambda_{\mu_1\mu_2})$.
A l'aide de la formule (2.11) on conclut facilement la relation suivante:
$$\left( \frac{\partial}{\partial t_0}\right)^{\ell} \tilde{\mathbf J}(\tilde t_0,t^{\prime}) = \Sigma^{(\ell)} \cdot \tilde{\mathbf J} (\tilde t_0, t^{\prime}). \leqno{(2.12)}$$
Ici
$\Sigma^{(\ell)}$
est une matrice
$M \times M$,
dont la  composante 
$\sigma^{(\ell)}_{ij},$ $\in End({\bf C}^{\mu})\otimes {\bf C}[t]$ $(i,j) \in [0,v_0]^2$
est d\'etermin\'ee comme suit:
$$\sigma^{(\ell)}_{ij} = \left\{ \begin{array}{cl}
\tilde S^{-\ell} \Lambda (\Lambda -id_{\mu}) \ldots (\Lambda-(\ell-1)id_{\mu}) & 0 \le i=j \le v_0 \\
i\ell \tilde S^{-(\ell +j-i)} \Lambda (\Lambda-id_{\mu}) \ldots (\Lambda-(\ell+j-i)id_{\mu}) & 0 \le j < i \le v_0 \\
0 & i < j \quad\mbox{o\`u}\quad \ell + j-i < 0 .\end{array} \right. \leqno{(2.13)}$$
C'est-\`a-dire:
$$\Sigma^{(\ell)} = \left( \begin{array}{llcl}
\sigma^{(\ell)}_{00} & \sigma^{(\ell)}_{01} & \cdots & \sigma^{(\ell)}_{0v_0} \\
\sigma^{(\ell)}_{10} & \sigma^{(\ell)}_{11} & \cdots & \sigma^{(\ell)}_{1v_0} \\
\vdots & & \ddots & \vdots \\
\sigma^{(\ell}_{v_00} & \sigma^{(\ell)}_{v_01} & \cdots & \sigma^{(\ell)}_{v_0v_0} \end{array} \right).$$
A l'aide de la matrice
$\Sigma^{(\ell)}$,
la condition (2.10) se transforme en l'\'equation suivante:
$$\langle \vec{r}, \Sigma^{(\ell)} \cdot \tilde{\mathbf J} (\tilde t_0, t^{\prime}) \rangle = 0 \leqno{(2.14)}$$
pour
$0 \le \ell \le \mu-1$. 
Puisque
$(\tilde t_0, t^{\prime}) \notin D$,
les valeurs propres de la matrice
$\tilde S (\tilde t_0,t^{\prime})$
sont sans multiplicit\'e. Par consequence,
$\tilde S (\tilde t_0, t^{\prime})$
est diagonalisable \`a l'aide d'une matrice $T$ \`a coefficients de fonctions semi-alg\'ebriques. C'est-\`a-dire:
$$T^{-1} \tilde S(\tilde t_0,t^{\prime})T = \textnormal{diag} (\tilde t_0-\tau_1(t^{\prime}), \tilde t_0 - \tau_2 (t^{\prime}), \ldots, \tilde t_0-\tau_{\mu}(t^{\prime})) \leqno{(2.15)}$$
o\`u
$\tau_i(t^{\prime}) \ne \tau_j(t^{\prime})$
pour
$i \ne j$.

En tenant compte de (2.12), (2.13), (2.14) on peut voir que l'existence d'une collection de nombres,d\'ependant
de $(\tilde t_0,t^{\prime},$ $\lambda_{00}, \cdots, \lambda_{\mu_1\mu_2}),$
$d^{(i)}_{M}(\tilde t_0,t^{\prime}),d^{(i)}_{M-1}(\tilde t_0,t^{\prime}), \ldots, d^{(i)}_0 
(\tilde t_0,t^{\prime}), i \ge 0$,
satisfaisant (2.16) ci-dessous est une condition 
suffisante pour que (2.10) pour
$\ell \ge \mu-1$
entra\^ine (2.10) pour
$\ell \ge \mu$
aussi:
$$d^{(i)}_M (\tilde t_0,t^{\prime}) \frac{{\tilde S}^{-(M-1)}}{(M+i)!} 
(\Lambda -(M+i-1)id_{\mu}) \ldots (\Lambda -i \cdot id_{\mu}) + 
\frac{d^{(i)}_{M-1}(\tilde t_0,t^{\prime})}{(M+i-1)!} {\tilde S}^{-M-i+1} (\Lambda -(\mu+i-2)id_{\mu}) \ldots \leqno{(2.16)}$$
$$\dots (\Lambda-i \cdot id_{\mu}) + \ldots + d^{(i)}_0 (\tilde t_0,t^{\prime}) \cdot \frac{\tilde S^{-i}}{i!} = 0,$$
avec $d^{(i)}_M (\tilde t_0,t^{\prime}) \ne 0$.
Vu (2.15), la relation (2.16)
est \'equivalente \`a l'\'equation suivante:
$$\sum^{M}_{\ell=1} d^{(0)}_{\ell} \frac{(\tilde t_0-\tau_{\nu}(t^{\prime}))^{-\ell}} {\ell !} \cdot \prod^{\ell-1}_{b=0} (\lambda_{\nu}-b) + d^{(0)}_0 = 0, \qquad \nu \in \mathbf I,$$
$$\sum^{M}_{\ell=0} \frac{d^{(i)}_{\ell}(\tilde t_0-\tau_{\nu}(t^{\prime}))^{-\ell}} {(\ell+i)!} \prod^{\ell+i-1}_{b=i} (\lambda_{\nu}-b) = 0$$
Celui-l\`a est, \`a son tour, \'equivalente \`a l'\'equation ci-dessous: 
$$\tilde{\Sigma}^{(i)} \cdot \vec{\mathbf d}^{(i)} = 0, \quad i \in \mathbf N$$
o\`u

$\tilde{\Sigma}^{(i)} = $
\begin{flushleft}
$
\left( \begin{array}{cccccccc}
{\scriptstyle\frac 1 {(M+i)!} \prod^{M+i-1}_{b=i}(\lambda_{00}-b),} & {\scriptstyle\frac{(\tilde t_0-\tau_1)} {(M+i-1)!} \prod^{M+i-2}_{b=i} (\lambda_{00}-b),} &\cdots, \\
{\scriptstyle \frac 1 {(M+i)!} \prod^{M+i-1}_{b=i} (\lambda_{01}-b),} & {\scriptstyle \frac{(\tilde t_0-\tau_2)} {(M+i-1)!} \prod^{M+i-2}_{b=i} (\lambda_{01}-b),} & \cdots,\\
\vdots & \vdots & \vdots \\
{\scriptstyle \frac {1} {(M+1)!} \prod^{M+i-1}_{b=i}(\lambda_{\mu_1\mu_2}-b)} 
& {\scriptstyle \frac{(\tilde t_0-\tau_{\mu})} {(M+i-1)!} 
\prod^{M+i-2}_{b=i} (\lambda_{\mu_1\mu_2}-b)} & 
\cdots,\\
{\scriptstyle \frac 1 {(M+i-1)!} \prod^{M+i-2}_{b=i} (\lambda_{00}-b)} & 
{\scriptstyle \frac{(\tilde t_0-\tau_1)} {(M+i-2)!} \prod^{M+i-3}_{b=i} 
(\lambda_{00}-b),} &\cdots  \\
\vdots & \vdots & \vdots \\ 
{\scriptstyle \frac 1 {(M+i-1)!} \prod^{M+i-2}_{b=i}(\lambda_{\mu_1\mu_2}-b)} 
& {\scriptstyle \frac{(\tilde t_0-\tau_{\mu})} {(M+i-2)!} \prod^{M+i-3}_{b=i} (\lambda_{\mu_1\mu_2}-b),} 
& \cdots,\\ 
\vdots & \vdots & \vdots\\ 
{\scriptstyle \frac 1 {(M+i-v_0)!} \prod^{M+i-v_0-1}_{b=i}(\lambda_{00}-b)} 
& {\scriptstyle \frac{(\tilde t_0-\tau_1)} {(M+i-v_0-1)!} \prod^{M+i-v_0-2}_{b=i} (\lambda_{00}-b),} 
& {\scriptstyle \cdots,}\\ 
\vdots & \vdots & \vdots \\ 
{\scriptstyle \frac 1 {(M+i-v_0)!} \prod^{M+i-v_0-1}_{b=i} (\lambda_{\mu_1\mu_2}-b)} 
& {\scriptstyle \frac{(\tilde t_0-\tau_{\mu})} {(M+i-v_0-1)!} \prod^{M+i-v_0-2}_{b=i} 
(\lambda_{\mu_1\mu_2}-b),} 
& {\scriptstyle \cdots,}\\ 
\end{array} \right. $
\end{flushleft}
\begin{flushright}
$
\left.   \begin{array}{cccccc}
\cdots, &\cdots,& \cdots,&\cdots & {\scriptstyle(\lambda_{00}-i) \frac{(\tilde t_0-\tau_1)^{M-1}} {(i+1)!},} & 
{\scriptstyle\frac{(\tilde t_0-\tau_1)^{M}} {i!}} \\
\cdots,& \cdots,& \cdots,&\cdots & {\scriptstyle (\lambda_{01}-i) \frac{(\tilde t_0-\tau_2)^{M-1}} {(i+1)!},} & 
{\scriptstyle \frac{(\tilde t_0-\tau_2)^{M}} {i!}} \\
\vdots & \vdots &\vdots &\vdots &\vdots
&\vdots \\
\cdots, &\cdots, &\cdots,&\cdots 
& {\scriptstyle (\lambda_{\mu_1,\mu_2}-i) \frac{(\tilde t_0-\tau_{\mu})^{M-1}} {(i+1)!},} 
& {\scriptstyle \frac{(\tilde t_0-\tau_\mu)^{M}} {i!}} \\
\cdots  &\cdots, &\cdots, &\cdots, & {\scriptstyle \frac{(\tilde t_0-\tau_1)^{M-1}} {i!},} & 0 \\
\vdots & \vdots & \vdots &\vdots & \vdots & \vdots \\
\cdots, &\cdots, &\cdots, &\cdots,& {\scriptstyle \frac{(\tilde t_0-\tau_{\mu})^{M-1}} {i!},} & 0 \\
\vdots & \vdots&&& & \vdots \\
{\scriptstyle \cdots,} & {\scriptstyle \frac{(\tilde t_0-\tau_1)^{M-v_0}} {i!},} &  0, &\cdots, &0, &0 \\
 \vdots & \vdots & \vdots& \vdots & \vdots & \vdots \\
{\scriptstyle \cdots,} & {\scriptstyle \frac{(\tilde t_0-\tau_{\mu})^{M-v_0}} {i!},} & 0 , &\cdots, &0,&0 
\end{array} 
\right). $
\end{flushright}
$$\vec{\mathbf d}^{(i)}= (d^{(i)}_{M}, d^{(i)}_{M-1}, \cdots, d^{(i)}_0) \in  \mathbf R^{M+1}, M = \mu (v_0+1).$$
Le vecteur 
$\vec{\mathbf d}^{(i)}$
est un normal aux $M$ vecteurs colonnes de la matrice
$\tilde{\Sigma}^{(i)}$.
Si on note ces $M$ vecteurs colonnes par
$\vec{v}^{(i)}_1,\cdots, \vec{v}^{(i)}_M$,
d'apr\`es alg\`ebre lin\'eaire on obtient l'expression suivante:
$$\vec{\mathbf d}^{(i)} = \textnormal{const.} \vec{v}^{(i)}_2 \wedge \cdots \wedge \vec{v}^{(i)}_M.$$
On en d\'eduit que
$$d^{(i)}_M (\tilde t_0(t'), \tau_1 (t'), \cdots, \tau_\mu (t'), \lambda_{00},
\cdots, \lambda_{\mu_1\mu_2}) 
= \tilde{\Sigma}^{(i)} \left( \begin{array}{l}
2,3, \cdots, M+1\\
1,2, \cdots, M \end{array} \right)$$
i.e. un
mineur $M \times M$ de la matrice
$\tilde{\Sigma}^{(i)}$.
Si on d\'efinit
$$\delta^{(i)}_M(t)= \prod_{\rho \in {\mathfrak S}_{\mu}, \sigma \in
{\mathfrak S}_{\bf I} }d^{(i)}_M 
( t_0, \tau_{\rho(1)}(t'), \cdots, \tau_{\rho(\mu)}(t'), \lambda_{\sigma(00)},
\cdots, \lambda_{\sigma(\mu_1\mu_2)}),$$
alors $\delta^{(i)}_M(t)$ est un polyn\^ome de degr\'e $\frac{M(M+1)}{2}
(\mu !)^2.$ Ici ${\mathfrak S}_{\bf I}$ indique le groupe de permutation de
$\mu$ indices de ${\bf I}.$ 

On note par $$B(t') = \prod_{1 \leq i <j\leq \mu}(\tau_i (t') -\tau_j (t'))^2 \leqno(2.17)$$
le polyn\^ome d\'efinissant l'ensemble de bifurcation. 
Alors il existent un entier $L \geq 0$ (eventuellement on calcule 
$L= \frac{(M \cdot(\mu-1)!^2)}{2}$) 
et des polyn\^omes $Q_{j}^{(i)}(\lambda_{00}, \cdots, 
\lambda_{\mu_1\mu_2}),$ $j=0, \cdots,$ 
$\frac{M(M+1)}{2\mu}(\mu !)^2 -L$ invariants sous l'action de 
${\mathfrak S}_{\bf I} $ tels que,
$$\delta^{(i)}_M(t) = 
(\sum_{\frac{(\mu -1)\ell}{2} +j =\frac{M(M+1)}{2\mu}(\mu !)^2-L } 
Q_j^{(i)}(\lambda)(\Delta_\mu(t))^j (B(t'))^\ell)\Delta_\mu(t)^L,
$$
o\`u $Q_0^{(i)}(\lambda) \not =0$ 
$Q_{\frac{M(M+1)}{2\mu}(\mu !)^2-L }^{(i)}(\lambda) \not =0. $ 

Evidemment 
$\delta^{(i)}_M(\tilde t_0,t^{\prime})=0$
pour les z\'eros
$(\tilde t_0, t^{\prime})$
de
$d^{(i)}_M(\tilde t_0, \tau_1, \cdots, \tau_\mu, \lambda_{00},
\cdots, \lambda_{\mu_1\mu_2})$
on a donc d\'emontr\'e l'\'enonc\'e du lemme  pour
$t = (t_0,t^{\prime})$
qui satisfait
$$\Delta_{\mu}(s(t)) \ne 0, \quad 
\delta^{(i)}_{M} (t) \ne 0, \quad i=0,1,2, \cdots.$$
Au voisinage de $t$ tel que	 $B(t')$ assez grand par rapport \`a $\Delta_{\mu}(s(t))$,  
il est un ouvert dense. 
Il est aussi un ouvert dense 
au voisinage de $t$ tel que $t_0$  assez grand par rapport \`a $B(t').$
En tenant compte de la quasihomog\'eneit\'e, c'est un emsemble ouvert dense de
${\bf R}^{\mu}.$ 
Cela d\'emontre le lemme.
C.Q.F.D.
                                                                                                                            \begin{remark} 
{\em Pour ${\bf J}(t)$ de (2.11),
B.Dubrovin \cite{Dub} a obtenu l'\'equation suivante:
$$ [\Delta_{\mu}(s(t))(\frac{\partial}{\partial t_0})^{\mu} 
+ g_{\mu-1}(t)(  \Lambda -id_\mu )
(\frac{\partial}{\partial t_0})^{\mu-1}+
g_{\mu-1}(t)(  \Lambda -id_\mu)(  \Lambda -2id_\mu)
\frac{\partial}{\partial t_0}^{\mu-2}
+\cdots $$
$$\cdots + g_1(t)(  \Lambda -id_\mu)\cdots (  \Lambda -(\mu-1)id_\mu )
\frac{\partial}{\partial t_0} +
(  \Lambda -id_\mu )\cdots (  \Lambda -\mu \cdot id_\mu )]
{\bf J}(t)=0,$$  pour la matrice diagonale $ \Lambda$
de (2.8) et des polyn\^omes $g_i(t)$ de degr\'e $i$ en $t_0$ t.q.
$$ det(\tilde S(t) + \lambda) = \Delta_\mu(s(t)) + \sum_{i=0}^{\mu-1}\lambda^{\mu-i}g_i(t).$$    
C'est une cons\'equence de l'\'equation (2.11) apr\`es application du th\'eor\`eme de Cayley-Hamilton
\`a la matrice $\tilde S(t).$
Cela implique que la multiplicit\'e des z\'eros d'une int\'egrale
 $$
\sum_{(\nu_1,\nu_2) \in \mathbf I} P_{\nu_1,\nu_2} J_{\nu_1\nu_2}(t)$$
avec $P_{\nu_1,\nu_2} \in \bf R$ ne doit pas d\'epasser $\mu$ hors de
l'ensemble critique $D =\{t; \Delta_{\mu}(s(t))=0\}.$

En fait, on peut d\'eduire de 
la d\'emonstration du Th\'or\`eme ~\ref{thm22} 
$iii)$ une estimation analogue pour $I_{x^ky^mdx}(s)$ 
. }
\label{remark21}
\end{remark}                                                                                                                                    
%%%%%%%%%%%%%%%%%%%%%REFERENCES%%%%%%%%%%%%%%%%%%%%%%%%%%%%%%%%%%

\vspace{\fill}

%%%%%%%%%%%%%%%%%%% ADDRESS %%%%%%%%%%%%%%%%

%

\noindent

\begin{flushleft}

 \begin{minipage}[t]{6.2cm}

  \begin{center}

{\footnotesize Moscow Independent University\\
Bol'shoj Vlasijevskij Pereulok 11,\\
  MOSCOW, 121002,\\

Russia\\

{\it E-mails}:  tanabe@mccme.ru, tanabesusumu@hotmail.com}

\end{center}

\end{minipage}\hfill

\end{flushleft}
\end{document}